\documentclass[a4paper,14pt]{article}

\usepackage{multicol}
\usepackage{amsthm}

\usepackage[dvips]{graphicx,color}
\usepackage{amssymb}
\usepackage{amsmath}

\begin{document}

\fontsize{14pt}{16.5pt}\selectfont

\begin{center}
\bf{Chaos on compact metric spaces \\
generated from symbolic dynamical systems}
\end{center}

\fontsize{12pt}{11pt}\selectfont
\begin{center}
Shousuke Ohmori\\ 
\end{center}

\noindent
\it{Faculty of Science and Engineering, Waseda~University, 3-4-1 Okubo, Shinjuku-ku, Tokyo 169-8555, Japan}\\
e-mail : 42261timemachine@ruri.waseda.jp\\
~~\\
\rm
\fontsize{11pt}{14pt}\selectfont\noindent

\noindent
{\bf Abstract}

We discuss Devaney chaos on compact metric spaces using a decomposition space characterized by topological nature of symbolic dynamics. A chaotic map obtained here is defined as a topologically conjugate of the chaotic map on a decomposition space which is induced by a chaotic map of symbolic dynamics. In particular, the chaotic character of the tent map and the baker map on $[0,1]$ are reconsidered based on decomposition dynamics involving symbolic dynamics with different two chaotic maps. As an example of compact metric space we exhibit a chaotic map existing on any given finite graph.

\section{Introduction}
\label{sec:1}

Chaotic dynamical systems have been widely studied in mathematics and physics and investigated its nature in a great deal of studies. Although there has been no generally mathematical or physical definition of chaos, the chaos defined by Devaney in his book\cite{Devaney} is widely known as presenting the essential feature of chaos. Symbolic dynamics has been used as one of the best tools to reveal various behaviors of dynamics with chaos such as Devaney chaos, consisting of a set $\{0,1\}^{\mathbb{N}}=\{\psi : \mathbb{N}\to \{0,1\}\}$ with a product 
topology $\tau _{0}^{\mathbb{N}}$ of $\tau _{0}=\{\phi,\{0\},\{1\},\{0,1\}\}$ and a continuous 
map $f$ from $\{0,1\}^{\mathbb{N}}$ into itself. Here we denote a symbolic dynamical system by a pair $(\{0,1\}^{\mathbb{N}},f)$. Note that the topological 
space $(\{0,1\}^{\mathbb{N}},\tau_{0}^{\mathbb{N}})$ is metrizable, namely, $\tau_{0}^{\mathbb{N}}=\tau_{d}$ where $d(\psi, \varphi)=\sum_{i=1}^{\infty}|\psi(i)-\varphi(i)|/2^i$ for $\psi,\varphi \in \{0,1\}^{\mathbb{N}}$.

The object of this article is to discuss Devaney chaos on compact metric spaces using a decomposition space of $(\{0,1\}^{\mathbb{N}},\tau_{0}^{\mathbb{N}})$. The decomposition space is exactly obtained as a usc decomposition by well-known fact (see {\bf 2-D} below) that to any compact metric space $Y$ there corresponds a decomposition space $\mathcal{D}$ of $(\{0,1\}^{\mathbb{N}},\tau_{0}^{\mathbb{N}})$ such that $\mathcal{D}$ is homeomorphic to $Y$. It will be shown that a chaotic map on a given compact metric space is defined as a topologically conjugate of a chaotic map on the decomposition space homeomorphic to the compact metric space and the chaotic map on the decomposition space is induced by the chaotic map of a symbolic dynamics satisfying a condition. In particular, the chaotic character of the tent map and the baker map on $[0,1]$ is reconsidered based on a decomposition space of symbolic dynamical systems with different two chaotic maps. As an example of a compact metric space on which the dynamics of the chaotic map can be concretely characterized, we exhibit a chaos existing for any given finite graph. 

\section{Definitions and Statements}
\label{sec:2}

Some definitions and statements which appear in the paper are listed below.

\begin{description}
\item[2-A.]\cite{Devaney}-\cite{Banks}
Let $(X,\tau)$ be a topological space and let $f:(X,\tau)\to (X,\tau)$ be a map. $f$ is said to be topologically transitive provided that for any pair of nonempty open sets $U,V$ of $X$, there exists $n>0$ such that $f^n(U)\cap V \not = \phi$, where $f^n$ denotes a composed mapping of $f$ with itself $n$ times ($f^0$ is the identity map). A orbit of $x\in X$ is a set $O^{+}(x)=\{x,f(x),f^2(x),\cdots\}$. $f$ has a dense orbit provided that $Cl O^+(x)=X$ for some point $x \in X$, where $Cl A$ stands for a closure of $A$, $A \subset X$. Note in a perfect T$_1$-space $(X,\tau)$ that if a map $f:(X,\tau)\to (X,\tau)$ has a dense orbit, then $f$ is topologically transitive. $f$ has a dense set of periodic points if $X=Cl Per(f)$ where $Per(f)=\{x\in X; f^n(x)=x,n>0\}$. Letting $f:X\to X$ and $g:Y\to Y$ be two maps, $f$ and $g$ are said to be topologically conjugate if there exists a homeomorphism $h:X \to Y$ such that the relation $g\circ h=h\circ f$ holds. If the map $h$ is continuous onto, they are said to be topologically semi-conjugate. Evidently, if $g$ is topologically semi-conjugate to $f$ i.e., $g\circ h=h\circ f$, $f$ having a dense set of periodic points and being topologically transitive, then $g$ is a topologically transitive and has a dense set of periodic points. Let $(X,\tau_d)$ be a metric space with a map $f:(X,\tau_d) \to (X,\tau_d)$. $f$ is said to be sensitive to initial conditions provided that there exists $\eta >0$ such that for any $x\in X$ and any neighbourhood $N(x)$ there is a point $y$ in $N(x)$ and $n\geq0$ satisfying $d(f^n(x),f^n(y))>\eta$. If an onto map $f:(X,\tau_d) \to (X,\tau_d)$ is topologically transitive, sensitive to initial conditions, and has a dense set of periodic points, then $f$ is said to be chaotic on $X$. However, if a continuous map $f$ is topologically transitive and has a dense set of periodic points, then $f$ turns out to be chaotic on an infinite metric space. It is noticed that continuity of $f$ can be replaced with the weaker condition in which $f$ is continuous at all points of an orbit of a periodic point. 
\item[2-B.]\cite{Nadler}
Let $(\mathcal{D},\tau(\mathcal{D}))$ be a decomposition space of a topological space $(X,\tau)$. That is, $\mathcal{D}$ is a set of nonempty, mutually disjoint subsets of $X$ such that $\bigcup \mathcal{D}=\bigcup _{D\in \mathcal {D}}=X$, and $\tau(\mathcal{D})$ is the decomposition topology defined by $\{\mathcal{U}\subset \mathcal{D}; \bigcup\mathcal{U}\in \tau \}$. A decomposition $\mathcal{D}$ of $(X,\tau)$ is said to be upper semi-continuous (usc) provided that whenever $D\in \mathcal{D}, U\in \tau$, and $D\subset U$, there is an open set $V$ with $D\subset V$ such that if $A\in \mathcal{D}$ and $A\cap V\not =\phi$, then $A \subset U$. We call a decomposition space $(\mathcal{D},\tau(\mathcal{D}))$ a usc decomposition space when $\mathcal D$ is usc. A subset $Y$ of $X$ is said to be $\mathcal{D}$-saturated if $Y$ is a union of subcollection of a decomposition $\mathcal{D}$, namely, $Y=\bigcup \mathcal{U}$ for some $\mathcal{U} \subset \mathcal{D}$. 

\noindent 
{\bf Statement}~{\it Let $(X,\tau)$ be a topological space, let $\mathcal{D}$ be a decomposition of $X$, and let $\pi : X \to \mathcal{D}$ be a natural map i.e., $\pi(x)=$ the unique $D\in \mathcal{D}$ such that $x \in \mathcal{D}$. Then, (i)-(iii) below are equivalent. 
\begin{description}
\item[(i)] $\mathcal{D}$ is a usc decomposition. 
\item[(ii)] If $D\in \mathcal{D}, U\in \tau,$ and $D\subset U$, there is an open set $V$ of $X$ such that $D\subset V \subset U$ and $V$ is $\mathcal{D}$-saturated. 
\item[(iii)] $\pi$ is a closed map from $(X,\tau)$ onto $(\mathcal{D},\tau(\mathcal{D}))$.
\end{description}
} 
{\bf Statement} {\it If $(X,\tau_d)$ is a compact metric space, then a usc decomposition space $(\mathcal{D},\tau(\mathcal{D}))$ of $(X,\tau_d)$ is compact and metrizable.}
\item[2-C.]\cite{Nadler,Engelking}Let $f : (X,\tau)\to (Y,\tau ')$ be a continuous closed onto map. Then,  the decomposition $\mathcal{D}_f = \{f^{-1}(y)\subset X;y\in Y\}$ is usc. Moreover, $h:(Y,\tau ') \to (\mathcal{D}_f,\tau(\mathcal{D}_f)), y\mapsto f^{-1}(y)$ is a homeomorphism. 
\item[2-D.]\cite{Nadler}Every nonempty compact metric space is a continuous image of a 0-dim, perfect, compact T$_2$-space.
\end{description}

\section{Map on a usc decomposition}
\label{sec:3}

Let us present at this section some results as to a map on a usc decomposition of a topological space, the map being induced by a map on the topological space. Throughout this section we assume that letting $G:X\to X$ be a map and letting $\mathcal{D}$ be a decomposition of $X$, the relation
\begin{center}
$(*)$~~$G(D)\in \mathcal{D}$ for any $D\in \mathcal{D}$
\end{center}
holds. We give the following lemmas which are easily verified from the definitions to begin with.\\

\noindent
{\bf Lemma 1.} {\it Let $G:X \to X$ be a map and let $\mathcal{D}$ be a decomposition (need not to be usc) of $X$. Define $H:\mathcal{D} \to \mathcal{D}, D\mapsto G(D)$. Then $(1)-(3)$ below hold.
\begin{description}
\item[(1)] If $G$ is onto, then $H$ is onto.
\item[(2)] If $G$ is one to one, then $H$ is one to one.
\item[(3)] If $\pi :X\to \mathcal{D}$ is a natural map, then the relation $H\circ \pi = \pi\circ G$ holds.
\end{description}}
~~\\
\noindent
{\bf Lemma 2.} {\it Let $G:(X,\tau) \to (X,\tau)$ be a map, let $\mathcal{D}$ be a usc decomposition space of $X$, and let $H:(\mathcal{D},\tau(\mathcal{D})) \to (\mathcal{D},\tau(\mathcal{D})), D\mapsto G(D)$. Then $(1)$ and $(2)$ below hold.
\begin{description}
\item[(1)] If $G$ is continuous, then $H$ is continuous.
\item[(2)] If $G$ is closed, then $H$ is closed.
\end{description}}

\begin{proof}
Let $D$ be an arbitrary point of $\mathcal{D}$ and let $\mathcal{U}\in \tau(\mathcal{D})$ with $H(D)\in \mathcal{U}$. By the definition of $\tau(\mathcal{D}), H(D) \subset \bigcup \mathcal{U}\in \tau$. Since $G$ is continuous and $(*)$ is satisfied, there exists $U \in \tau$ with $D\subset U$ such that whenever $x \in U$, $G(x) \in \bigcup \mathcal{U}$. Since $\mathcal {D}$ is usc, there is an open set $V$ of $X$ such that $D\subset V \subset U$ and $V$ is $\mathcal{D}$-saturated. That is, there is a subset $\mathcal{V}$ of $\mathcal{D}$ such that $V=\bigcup \mathcal{V}$. Then, $D\in \mathcal{V} \in \tau(\mathcal{D})$ and it follows that $D'\in \mathcal{V}$ implies $H(D')=G(D')\in \mathcal{U}$. Thus, $H$ is continuous. Next, let $\mathcal{K}\in \Im(\mathcal{D})$. Then, $\bigcup \mathcal{K}\in \Im$. Since $G$ and the natural map $\pi$ are closed from {\bf 2-B.}, it is easily verified from (3) of Lemma 1 that $H(\mathcal{K})=\pi(G(\bigcup \mathcal{K}))\in \Im(\mathcal{D})$. Therefore, $H$ is closed.
\end{proof}
~~\\
\noindent
{\bf Corollary 3.} {\it If $G:(X,\tau) \to (X,\tau)$ is homeomorphic and $\mathcal{D}$ is a usc decomposition space of $X$, then $H:(\mathcal{D},\tau(\mathcal{D}))\to (\mathcal{D},\tau(\mathcal{D})), D\mapsto G(D)$ is homeomorphic.}
~~\\

Define $H$ as in Lemma 2. By (3) of Lemma 1, $G$ and $H$ are topologically semi-conjugate with respect to the natural map $\pi:X\to \mathcal{D}$. Therefore, it immediately follows that\cite{Banks1997,Balibrea}\\

\noindent
{\bf Lemma 4.} {\it If $G:(X,\tau) \to (X,\tau)$ is topologically transitive and has a dense set of periodic points, then $H:(\mathcal{D},\tau(\mathcal{D})) \to (\mathcal{D},\tau(\mathcal{D})), D\mapsto G(D)$is topologically transitive and has a dense set of periodic points.}
~~\\

\noindent
{\bf Corollary 5.} {\it Let $X$ be a compact metric space. If $G:X \to X$ is chaotic on $X$ and $\mathcal{D}$ is an infinite usc decomposition space of $X$, then $H:\mathcal{D} \to \mathcal{D}, D\mapsto G(D)$ is chaotic on $\mathcal{D}$.}
\begin{proof}
The topologically transitivity and the existence of a dense set of periodic points of $H$ are shown from Lemma 4. Since $X$ is a compact metric space, from {\bf 2-B.}, the usc decomposition space $\mathcal{D}$ is a compact metric space. Therefore, it is easy to verified that $H$ is chaotic on $\mathcal{D}$.  
\end{proof}

\section{Chaos on compact metric spaces generated from symbolic dynamical systems}
\label{sec:4}

\subsection{General Consideration}
\label{sec:4.1}

Let $\{0,1\}^{\mathbb{N}}(=(\{0,1\}^{\mathbb{N}},\tau_{0}^{\mathbb{N}}))$ be a topological space as stated in \S\ref{sec:1}. Since it is easily proved that $\{0,1\}^{\mathbb{N}}$ is a 0-dim, perfect, and compact metric space, for any compact metric space denoted by $(Y,\tau_\rho)$ there exists a continuous onto map $f:(\{0,1\}^{\mathbb{N}},\tau_{0}^{\mathbb{N}}) \to (Y,\tau_\rho)$ ({\bf 2-D.}). Since $\{0,1\}^{\mathbb{N}}$ is compact and $Y$ is T$_2$, $f$ is closed. Thus, by {\bf 2-C.}, there exists a usc decomposition space $(\mathcal{D}_f,\tau(\mathcal{D}_f))$ of $\{0,1\}^{\mathbb{N}}$ such that $h:(Y,\tau _\rho) \to (\mathcal{D}_f,\tau(\mathcal{D}_f)), y\mapsto f^{-1}(y)$ is a homeomorphism. Assume $G: \{0,1\}^{\mathbb{N}} \to \{0,1\}^{\mathbb{N}}$ to be a continuous onto map satisfying $(*)$ in \S\ref{sec:3}. From Lemma 1 and 2, then, $H:\mathcal{D} \to \mathcal{D},D\mapsto G(D)$ is a continuous onto map. Hence, the continuous onto map $F:(Y,\tau_\rho)\to (Y,\tau_\rho), y\mapsto h^{-1} \circ H \circ h(y)$ is obtained, which is topologically conjugate to $H$. That is, we can derive a discrete dynamical system composed of a pair of $(Y,F)$ from the decomposition dynamics of $(\mathcal{D},H)$. In particular, if $G$ is chaotic on $\{0,1\}^{\mathbb{N}}$, then $H$ is chaotic by Corollary 5. Therefore, a chaotic map $F$ on an infinite compact metric space $Y$ is obtained systematically from the symbolic dynamical system ($\{0,1\}^{\mathbb{N}},G$). 

However, it is generally difficult to specify a continuous map $G$ on $\{0,1\}^{\mathbb{N}}$ as it holds $(*)$, for the condition $(*)$ depends on the usc decomposition $\mathcal{D}$ of $\{0,1\}^{\mathbb{N}}$ which is homeomorphic to a given compact metric space. As one of the methods to avoid such difficulty, we provide the following lemma regarding the topologically transitivity and existence of a dense set of periodic points.\\
    
{\bf Lemma 6.} {\it Let $G:(X,\tau) \to (X,\tau)$ be a map and let $\mathcal{D}$ be a decomposition space having a subset $\mathcal{D}'$ of $\mathcal{D}$ such that $G(D) \not \in \mathcal{D}$ if and only if $D\in \mathcal{D}'$. Let $H:(\mathcal{D},\tau(\mathcal{D}))\to (\mathcal{D},\tau(\mathcal{D}))$ be a map defined by
\begin{center}~~~
$H(D)=\left\{
\begin{array}{lcr}
G(D)~,~~D\in \mathcal{D}-\mathcal{D}'\\
D'~~~~~,~~D\in \mathcal{D}'\end{array}
\right.$
\end{center} 
where $D'$ is chosen as a point of $\mathcal{D}$ for which $H$ is onto. Then $(1)$ and $(2)$ below hold.
\begin{description}
\item[(1)] If $G$ has a dense set of periodic points and $H\circ \pi(\psi) = \pi\circ G(\psi)$ holds for each $\psi \in Per(f)$, then $H$ has a dense set of periodic points.
\item[(2)] If $G$ has a dense orbit $O^+(x_0)$ and $H\circ \pi(\psi) = \pi\circ G(\psi)$ holds for each $\psi \in O^+(x_0)$, then $H$ has a dense orbit, in fact, $\mathcal{D}=Cl O^+(\pi(x_0))$.
\end{description}}
\begin{proof}
We will show (1). Let $D\in \mathcal {D}$ and let $\mathcal{U}\in \tau(\mathcal{D})$ with $D\in \mathcal{U}$. By definition of $\tau(\mathcal{D})$, $D\subset \bigcup \mathcal{U}\in \tau$. Let $\psi \in D$. By $ClPer(G)=X$ there exist $n>0$ and $\psi_0 \in \bigcup \mathcal{U}$ such that $G^n(\psi_0)=\psi_0$. By assumption of (1), $H^n(\pi(\psi_0))=\pi(\psi_0)$ and hence $\pi(\psi_0)\in \mathcal{U}\cap Per(H)$. In the same way, (2) can be proven. Therefore, the proof of Lemma 6 is complete.
\end{proof}

According to practical construction of a chaotic map on a compact metric space we consider, it is necessary to clarify the concrete description of the usc decomposition which is homeomorphic to the compact metric space. We have investigated the definite form of use decompositions of $\{0,1\}^{\mathbb{N}}$ involved with compact metric spaces consisting of some arcs\cite{Ohmori2017}. Therefore, in the following subsection we confine compact metric spaces to this type.


\subsection{Reconsideration of chaos on [0,1]}
\label{sec:4.2}

As a simple example of a compact metric space consisting of arcs, the closed interval $[0,1]$ is considered. It is familiar with the two chaotic maps on $[0,1]$, i.e., the tent map and the baker map. First we focus on the tent map that is defined by $T:[0,1] \to [0,1], T(x)=2x$ if $0\leq x \leq 1/2$, and $T(x)=2(1-x)$ if $1/2\leq x \leq 1$. By {\bf 2-C.} and {\bf 2-D.} there exists a decomposition $\mathcal{D}$ of $\{0,1\}^{\mathbb{N}}$ such that $h:[0,1] \to \mathcal{D},y\mapsto f^{-1}(y)$ where $f$ is a continuous map from $\{0,1\}^{\mathbb{N}}$ onto $[0,1]$. Here $\mathcal{D}$ can be taken as satisfying that for any $y\in [0,1]$ there uniquely corresponds to $X_{k_1k_2\cdots}\in \mathcal{D}$ such that $f^{-1}(y)=X_{k_1k_2\cdots}$ where $k_1,k_2,\cdots \in \{0,1\}$ are values occurring when $y$ is represented by using a binary system, namely, $y=\sum_{i=1}^{\infty}k_i/2^i$, and subspaces $\{\psi \in \{0,1\}^{\mathbb{N}};\psi(i)=k_i,i=1,\cdots,n\}$, $\{\psi \in \{0,1\}^{\mathbb{N}};\psi(i)=k_i,i\in \mathbb{N}\}$ of $\{0,1\}^{\mathbb{N}}$ are denoted by $X_{k_1k_2\cdots k_n}$, $X_{k_1k_2\cdots}$, respectively. Note that if a point $y$ of $[0,1]$ is contained in $M=\{l/2^n\in [0,1];n=1,2,\cdots, l=1,2,\cdots ,2^n-1\}$, then there are finitely many points $k_1,\cdots,k_n$ of $\{0,1\}$ such that $f^{-1}(y)=X_{k_1k_2\cdots k_n100\cdots}\cup X_{k_1k_2\cdots k_n011\cdots}$. To induce the tent map, we adopt as a map of symbolic dynamics $C:\{0,1\}^{\mathbb{N}}\to \{0,1\}^{\mathbb{N}} $ defined by  $(C(\psi))(i)=\psi(i+1)$ for any $i$ if $\psi(1)=0$ and $(C(\psi))(i)=1-\psi(i+1)$ for any $i$ if $\psi(1)=1$. It follows that $C$ is a continuous onto map not to be one to one. Note that $C$ satisfies $(*)$ in \S\ref{sec:3}. In fact, let $D\in \mathcal{D}_f$. Then, there is a point $y\in [0,1]$ such that $f^{-1}(y)=D$. Assume $y \not \in M$. There are infinitely many elements $k_1,k_2,\cdots$ of $\{0,1\}$ such that $y=\Sigma _{i=1}^{ 	\infty} k_i/2^i$. If $k_1=0$, then $C(D)=C(X_{k_1k_2k_3\cdots})=X_{k_2k_3\cdots}$. If $k_1=1$, then $C(D)=X_{1-k_21-k_3\cdots}$. In the both cases, $C(D)\in \mathcal{D}$. In the same way, $C(D) \in \mathcal{D}$ is proved for $y\in M$. Therefore, it is possible to define $H:\mathcal{D}_f \to \mathcal{D}_f$ by $H(D)=C(D)$ and then $H$ is a continuous onto map by Lemma 1, 2. Since it is easily shown that $C$ is chaotic on $\{0,1\}^{\mathbb{N}}$, by \S\ref{sec:4.1}, the chaotic map $F:[0,1]\to [0,1]$ is obtained where $F(y)=h^{-1} \circ H \circ h(y)$. Note that the map $F$ is equal to the tent map $T$ and hence we lead the common knowledge that the tent map is chaotic on $[0,1]$.

Next we consider the baker map $B:[0,1]\to [0,1]$ which is defined by $B(x)=2x$ if $0\leq x \leq 1/2$ and $B(x)=2x-1$ if $1/2< x \leq 1$. Let $S$ denote a shift map of $\{0,1\}^{\mathbb{N}}$, i.e., $(S(\psi))(i)=\psi(i+1)$ for any $i\in \mathbb{N}$. Then, $S$ is a continuous onto and chaotic map. Note $S$ does not satisfy $(*)$ in \S\ref{sec:3} at the point $f^{-1}(1/2)$. Taking lemma 6 into consideration, define a map $H':\mathcal{D} \to \mathcal{D}$ by $H'(D)=S(D)$ for each $D\in \mathcal{D}-\{ f^{-1}(1/2) \}$ and $H'(f^{-1}(1/2))=X_{111\cdots}=f^{-1}(1)$. Since $Per (S)\cap f^{-1}(1/2)=O^+(\psi_0)\cap f^{-1}(1/2)=\phi$ and $H'\circ \pi(\psi)=\pi \circ S(\psi)$ for any $\psi \in Per(S)\cup O^+(\psi_0)$ where $\psi_0 = (01;00011011;000\dots)$, lemma 6 can be used to show that $H'$ is chaotic on $\mathcal{D}$. Obviously, $H'$ is not continuous at $f^{-1}(1/2)\in \mathcal{D}$. Since $B=h^{-1}\circ H'\circ h$ is clearly equal to the baker map, we derive the common knowledge that the baker map is chaotic on $[0,1]$ from the symbolic dynamical system $(\{0,1\}^{\mathbb{N}},S)$. It is noticed that $S$ and $C$ are topologically conjugate. In fact, let $R:\{0,1\}^{\mathbb{N}}\to \{0,1\}^{\mathbb{N}}$ be a map defined by $R(\psi)(i)=C^i(\psi)(1)$ for $i\in \mathbb{N}$. Then, it is clear that $R$ is homeomorphism satisfying a relation $S\circ R=R\circ B$. 

\subsection{Chaos on a finite graph}
\label{sec:4.3}

The discussion of \S\ref{sec:4.1} can be available to the topological space constructing of some arcs such as a finite graph. In the subsection we devote to create a chaotic map on a given finite graph. A finite graph $(E,\tau_ \rho )$ is defined from topology as a continuum that can be written as the union of finitely many arcs $E_1,\cdots,E_r$ any two of which are either disjoint or intersect only in one or both of their end points. Since each $E_i$ is an arc with end points $a_i$ and $b_i$, there is a homeomorphic map $h_i$ from $[0,1]$ onto $E_i$ such that $h_i(0)=a_i$ and $h_i(1)=b_i$. We assume without loss of generality that each end point is a node of the graph and the edge between two nodes is an open arc. That is, $N(E)=\{a_i,b_i;i=1,\dots,r\}$ and $e(E)=\{E_1-\{a_1,b_1\},\cdots,E_r-\{a_r,b_r\}\}$ are finite sets of nodes and edges of $E$, respectively. 

Now let us consider a symbolic dynamical system $\{0,1\}^{\mathbb{N}}$ with shift map $S$. By using {\bf 2-D} it is confirmed that there is a partition $\{ X_{0},X_{10},\cdots,X_{11\cdots^{r-1}10},X_{11\cdots^{r-1}11}\}$ of $\{0,1\}^{\mathbb{N}}$ such that each arc $E_i$ is a continuous image of $X_{11\cdots ^i10}$ where $X_{11\cdots ^i10}$ stands for $\{\psi \in \{0,1\}^{\mathbb{N}};\psi(1)=1,\dots,\psi(i-1)=1,\psi(i)=0 \}$. Then we can obtain by {\bf 2-C} a usc decomposition $\mathcal{D}$ of $\{0,1\}^{\mathbb{N}}$ which is homeomorphic to $E$ by a homeomorphism $h$ and satisfies the followings; There correspond subsets $\mathcal{D}_1,\cdots, \mathcal{D}_r$ of $\mathcal{D}$ to $E_1,\cdots, E_r$ such that $\mathcal{D}=\bigcup_{i=1}^r\mathcal{D}_i$, $|\mathcal{D}_i\cap \mathcal{D}_j|\leq 1$ for $i\not =j$, and for each $i$ and for any $y\in E_i$ there is a unique point $X_{11\cdots^i10;k_1k_2\cdots}$ of $\mathcal {D}_i$ such that $h(y)=X_{11\cdots^i10;k_1k_2\cdots}$ where $k_1,k_2,\cdots$ are obtained when $h_i^{-1}(y)\in [0,1]$ is represented by using a binary system and $X_{11\cdots^i10;k_1k_2\cdots}=\{\psi \in \{0,1\}^{\mathbb{N}};\psi(1)=1,\dots,\psi(i-1)=1,\psi(i)=0,\psi(i+j)=k_j$ for $j\in \mathbb{N} \}$.
 It is found some points of $\mathcal{D}$ at which $S$ cannot satisfy the condition $(*)$ in \S\ref{sec:3} as well as the case of the baker map we have shown in \S\ref{sec:4.2}. For instance, let $D'_r$ be a point of $\mathcal{D}_r$ such that $h^{-1}_r(h^{-1}(D'_r))=1/2$, namely, $D'_r=X_{11\cdots^{r-1}11;100\cdots}\cup X_{11\cdots^{r-1}11;011\cdots}$. Then, $S(D'_r)=X_{11\cdots^{r-1}11;00\cdots}\cup X_{11\cdots^{r-1}10;11\cdots} \not \in \mathcal{D}$. Note that the same difficulty occurs if we adopt as a continuous map of symbolic dynamics a map $C$ in \S\ref{sec:4.2} in place of $S$. To add to $D'_r$ of $\mathcal{D}_r$, $S(D'_1)\not \in \mathcal{D}$ is easily shown for the point $D'_1$ of $\mathcal{D}_1$ being $h^{-1}_1(h^{-1}(D'_1))=1/2$. Furthermore, it is not necessarily for $S$ to satisfy $(*)$ at each point of $\mathcal{D}$ corresponding to a node of the graph $E$, a node being a point which is regarded as an end point of some arcs in $\{E_j;j=1,\cdots,r\}$. For instance, a point $D'$ of $\mathcal{D}$ corresponding a node $y(=h^{-1}(D))\in E$ with just three edges $E_l-\{y,e_l\}$, $E_m-\{y,e_m\}$, and $E_n-\{y,e_n\}$ ($l\not =m\not =n, l,m,n\in \{ 2,\cdots, r-1\}$) is $X_{11\cdots^{l}10;00\cdots}\cup X_{11\cdots^{m}10;11\cdots}\cup X_{11\cdots^{n}10;00\cdots}$ where $h_l(0)=h_m(1)=h_n(0)=y, h_l(1)=e_l, h_m(0)=e_m, h_n(1)=e_n$. Then, $S(D)=X_{11\cdots^{l-1}10;00\cdots}\cup X_{11\cdots^{m-1}10;11\cdots}\cup X_{11\cdots^{n-1}10;00\cdots}$ is not necessarily to be contained in $\mathcal{D}$. Letting $\mathcal{D}'=\{D'_1,D'_r, h(y);y\in N(E)\}$, it is verified that$S(D)\in \mathcal{D}$ if and only if $D\in \mathcal{D}- \mathcal{D}'$. Thus, define $H:\mathcal{D}\to \mathcal{D}$ by $H(D)=S(D)$ for any $D\in \mathcal{D}- \mathcal{D}'$ and $H(D)=D$ for any $D\in \mathcal{D}'$. Since $H\circ \pi(\psi) = \pi\circ S(\psi)$ holds for each $\psi \in Per(f)\cup O^+(\psi_0)$, from Lemma 6, $H$ has a dense set of periodic points and has a dense orbit. Therefore, $F:E\to E,y\mapsto h^{-1} \circ H \circ h(y)$ is a chaotic map on $E$. 

Finally, one describes briefly orbits of $F$ on $E$. By the definition of $H$, each node and the two points $y'_1,y'_r$ corresponding to $h_1(y'_1)=h_r(y'_r)=1/2$ are fixed points of $E$ for $F$. Let $p_i$ be a point in an arc $E_i,i\not=1,r$ such that $h_i(p_i)=t\in [0,1]-\{0,1\}$. Then, $F$ transforms $p_i$ into $p_{i-1}$ which is contained in $E_{i-1}$ such that $h_{i-1}(p_{i-1})=t$. Acting $F$ to $p_i$ successively, a sequence $\{p_i, p_{i-1},\dots,p_1\}$ in $E$, each point $p_j$ of which is contained in $E_i$ holding $h_j(p_j)=t$ up to $j=1$, is obtained. $F$ transforms again $p_1$ in $E_1$ into a point $q$ in $E$ such that if $q$ is in $E_s$ and $t=\sum_{i=1}^\infty k_i/2^i$ where $k_1,k_2,\cdots \in \{0,1\}$, then $h_s(q)=\sum_{i=1}^\infty k_{i+1}/2^i$. As a whole $F$ shows chaotic orbit on $E$.

\section{Conclusion}
We discuss Devaney chaos on compact metric spaces using a usc decomposition space $\mathcal{D}$ obtained by topological nature of symbolic dynamical systems $(\{0,1\}^{\mathbb{N}},f)$. The chaotic map $F$ on a given infinite compact metric space $Y$ can be derived as a topologically conjugate of a chaotic map $H$ on $\mathcal{D}$ with the homeomorphism $h$ from $Y$ onto $\mathcal{D}$. Here since $\mathcal{D}$ is usc, using lemmas in \S\ref{sec:3}, $H$ is induced by the chaotic map $f$ of a symbolic dynamics if $f$ satisfying the condition $(*)$. The flow of construction of $F$ is summarized in \S\ref{sec:4.1}. In particular, the tent map and the baker map on $[0,1]$ are obtained based on symbolic dynamical systems $(\{0,1\}^{\mathbb{N}},C)$ and $(\{0,1\}^{\mathbb{N}},S)$, respectively. According to the case of the baker map, it is necessary to apply lemma 6 to $(\{0,1\}^{\mathbb{N}},S)$ because of existence of the points at which $(*)$ is violated for $S$. As an example of a compact metric space which consists of arcs, we create a chaotic map $F$ existing on any given finite graph generated from a symbolic dynamics $(\{0,1\}^{\mathbb{N}},S)$. Since the existence of $F$ is independent of the connectedness of graph, it must be shown the chaos on a compact metric space such as a disjoint union of some arcs. It is emphasized again that the construction of a chaotic map stated in \S\ref{sec:4.1} can be applied to any compact metric space.

\bigskip

\noindent
{\bf Acknowledgment}\\
The author is grateful to Prof. T. Yamamoto, Prof. Y. Yamazaki, and Prof. Emeritus A.Kitada at Waseda university for useful suggestions and encouragements.


\begin{thebibliography}{9}
\bibitem{Devaney} R. Devaney, \emph{An Introduction to Chaotic Dynamical Systems, 2ed edition} (Addison-Wesley 1989).
\bibitem{Holmgren} R. A. Holmgren, \emph{A First Course in Discrete Dynamical Systems, 2ed edition} (Springer 2000).
\bibitem{Silverman} S. Silverman, Rocky. Mt. J. Math. {\bf 22} 353 (1992).
\bibitem{Banks} J. Banks, J. Brooks, G. Cairns, G. Davis, and P. Stacey, Amer. Math. Mon. {\bf 99} 332 (1992).
\bibitem{Nadler} S. B. Nadler, Jr. \emph{Continuum theory.} (Marcel Dekker, 1992). 
\bibitem{Engelking} R. Engelking, \emph{General Topology} (Heldermann Verlag Berlin 1989).
\bibitem{Banks1997} J. Banks, Ergod. Th. \& Dynam. Sys. {\bf 17} (1997) 505.
\bibitem{Balibrea} F. Balibrea, T. Downarowicz, R. Hric, L'. Snoha, and V. Spitalsky, Ergod. Th. \& Dynam. Sys. {\bf 29} (2009) 737.
\bibitem{Ohmori2017} S. Ohmori, T. Yamamoto, and A. Kitada, arXiv:1708.02748v1 [math-ph] 9 Aug 2017.
\end{thebibliography}
\end{document}